\documentclass[]{article}
\usepackage{latexsym,epsf,rotating}
\usepackage[dvips]{epsfig}

\date{ }
\begin{document}

\title{Characterization of the Distribution of Twin Primes}
\author{P.F.~Kelly\footnote{patrick\_kelly@ndsu.nodak.edu}
and Terry~Pilling\footnote{terry@mailaps.org} \\
Department of Physics \\
North Dakota State University \\
Fargo, ND, 58105-5566 \\
U.S.A.}
\maketitle

\begin{abstract}
We adopt an empirical approach to the characterization of
the distribution of twin primes within the set of primes, 
rather than in the set of all natural numbers.
The occurrences of twin primes in any finite sequence of primes 
are like fixed probability random events.
As the sequence of primes grows, the probability decreases 
as the reciprocal of the count of primes to that point.
The manner of the decrease is consistent with the Hardy--Littlewood 
Conjecture, the Prime Number Theorem, and the Twin Prime Conjecture.
Furthermore, our {\it probabilistic model}, is simply parameterized.
We discuss a simple test which indicates the consistency of the
model extrapolated outside of the range in which it was constructed.
\end{abstract}

\noindent
Key words: Twin primes

\noindent
MSC: 11A41 (Primary), 11Y11 (Secondary)

\section{Introduction}

Prime numbers~\cite{cohen}, with 
their many wonderful properties, have been an intriguing subject 
of mathematical investigation since ancient times.
The ``twin primes,'' pairs of prime numbers $\{ p, p+2 \}$ are a 
subset of the primes and themselves possess remarkable properties.
In particular, we note that the {\it Twin Prime Conjecture},
that there exists an infinite number of these prime number pairs
which differ by 2, is not yet proven~\cite{hardy,marzinske}.

In recent years much human labor and computational effort have been 
expended on the subject of twin primes.
The general aims of these researches have been three-fold:
the task of enumerating the twin primes~\cite{halberstam} 
({\it i.e.,} identifying the members of this particular subset of 
the natural numbers, and its higher-order variants ``$k$-tuples'' of primes), 
the attempt to elucidate how twin primes are distributed among the 
natural numbers~\cite{huxley,brent1,brent2,wolf1} ({\it especially} searches 
for long gaps in the sequence~\cite{nicely1,wolf2,odlyzko}),
and finally, the precise estimation of the value of Brun's 
Constant~\cite{wolf3}.

Many authors have observed that the twin primes, along with the primes
themselves, generally become more sparse or diffuse as their magnitude
increases.
In fact the {\it Prime Number Theorem} may be rephrased to state that 
the number of prime numbers less than or equal to some large 
(not necessarily prime) number $N$ is 
approximately\footnote{In fact, more accurate 
approximations are known, but the formulae we quote suffice for 
our purposes.}
\begin{equation}
\pi_{1}(N) \sim \int_2^N \frac{1}{\ln(x)} dx\, .
\label{pi1}
\end{equation}
A similar result is believed to hold for the number of twin primes 
where each element of the pair is less than or equal to large $N$
\begin{equation}
\pi_{2}(N) \sim 2 c_2 \int_2^N \frac{1}{\left(\ln(x) \right)^2} dx\, ,
\label{pi2}
\end{equation}
where the ``twin prime'' constant~\cite{tpconst1} $c_2$ has the 
numerical value $c_2 = 0.661618\ldots$, and is currently known to 
many decimal places~\cite{tpconst2,tpconst3}.
The expression~(\ref{pi2}) is the first instance of the Hardy--Littlewood
conjectures which estimate the multiplicity of $k$-tuples of primes 
smaller than natural number $N$~\cite{hardy}.

In our investigation, we sought to account for the effect of the 
primes themselves becoming more rarefied by examining the distribution
of twin primes within the prime numbers.
We have done this for the set of prime numbers less than 
approximately $4 \times 10^9$, and have observed that within this range
the occurrence of twin primes may be characterized as  
{\it slowly-varying-probability random} events.

\section{Method and Results}  

To ensure that there is no confusion, let us first make very clear 
our methodology.
We generated prime numbers in sequence, {\it viz,}
$\ P_1 = 2, P_2 = 3, P_3 = 5 \ldots $,
and within this sequence identified twin primes and their
{\it prime separations} as illustrated
somewhat schematically below.
\[
\cdots \ P_i \ \left( P_{i+1} \ P_{i+2} \right) \ P_{i+3} \ P_{i+4} 
\ \left( P_{i+5} \ P_{i+6} \right) \ P_{i+7} \ P_{i+8} \ P_{i+9} \ P_{i+10}
\ \left( P_{i+11} \ P_{i+12} \right) \ P_{i+13} \ \cdots  
\]
We say that the first pair of twins in the above sequence has a 
prime separation of 2.
There are two non-twin prime numbers, {\it i.e., singletons}, which 
occur between the second prime element, $P_{i+2}$, of the first 
twin and the first prime element, $P_{i+5}$, of the subsequent twin.
Similarly, the second pair of twins has prime separation 4.
Note that there are many twins with prime separation equal to zero:
for example 
$(5\ 7) (11\ 13)$, or $(137\ 139) (149\ 151)$.
Actually, all of the prime 4-tuples $\left( P, P+2, P+6, P+8 \right)$
are comprised of a pair of twins with zero separation in primes.

There is an irregularity with our definition of separation
for the first few primes: $2\  (3\ 5) (5\ 7)$,
where the pairs in fact overlap, yielding a prime separation of $-1$.
Fortunately, for very well-known reasons such overlapping twins
do not ever recur and we choose to begin our analysis with the twin 
$(5\ 7)$.
For instance, in the set of seven twins between 5 and 100,
\[
\left\{ (5\ 7) (11\ 13) (17\ 19) (29\ 31) (41\ 43) (59\ 61) 
(71\ 73) \right\},
\]
there are 6 separations.
Three of these happen to be 0, two are 1, and one is 2, so the 
relative frequencies for separations $s = 0 , 1 , 2$ are
$\frac{1}{2} , \frac{1}{3}$, and $\frac{1}{6}$, respectively.

\begin{figure}[htb]
\begin{center}
\begin{turn}{-90}
\leavevmode
\epsfxsize=3.5in
\epsfysize=5.5in
\epsfbox{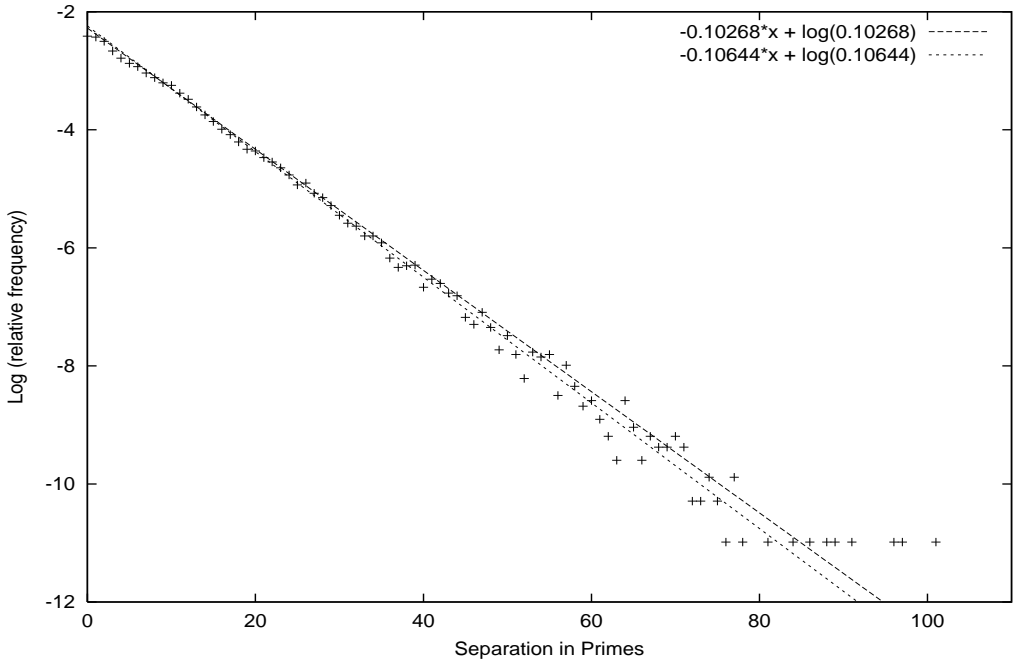}
\end{turn}
\end{center}
\caption{Data and Linear Fits for log(frequency) {\it vs.}
separation, in the case $N = 10^7$.}
\label{graph1}
\end{figure}

From the set of primes, all separations between pairs of twins 
up to a fixed number $N$ were computed and tabulated.
[We chose certain values of $N$ in the range $79561$ to $4020634603$.
These particular numbers are the second prime elements of the 
thousandth and twelve-millionth twins respectively.
Many of our $N$'s were chosen such that our analysis started and 
ended on twins.
However the behavior that we observe holds for all $N$, sufficiently
large, with the understanding that the singleton primes between the 
last twin and the upper bound $N$ are ignored.]
The logarithm of the relative frequency of occurrence of each 
separation in each of our 
analyses appears to obey a surprisingly simple relation
as illustrated {\it schematically} in Figure~\ref{graph1}.

Two comments must be made.
The first is that this remarkable behavior is perfectly characteristic
of a completely {\it random} system.
We infer that as one approaches {\it each} prime member in the 
sequence of primes following a twin, the likelihood of it being the first 
member of the next twin prime is {\it constant}!
By way of analogy, we can consider a radioactive substance.
The likelihood of one of its atoms decaying in any short time 
interval is fixed, with the effect that the probability that 
the next decay occurs at time $t$ is just 
${\cal N} e^{-\gamma t}$, where $\gamma$ is the decay rate, and ${\cal N}$
is a constant to ensure appropriate normalization.
The measured slope of the line fit to our data provides a decay constant
which is particular to the twin primes.
Again, recourse to our analogy is warranted.
The decay rate of a radioactive element is a {\it defining characteristic} 
of that element.
We appear to have the occurrences of twin primes governed by 
a similar sort of ``prime constant.''

\begin{figure}[htb]
\begin{center}
\begin{turn}{-90}
\leavevmode
\epsfxsize=3.5in
\epsfysize=5.5in
\epsfbox{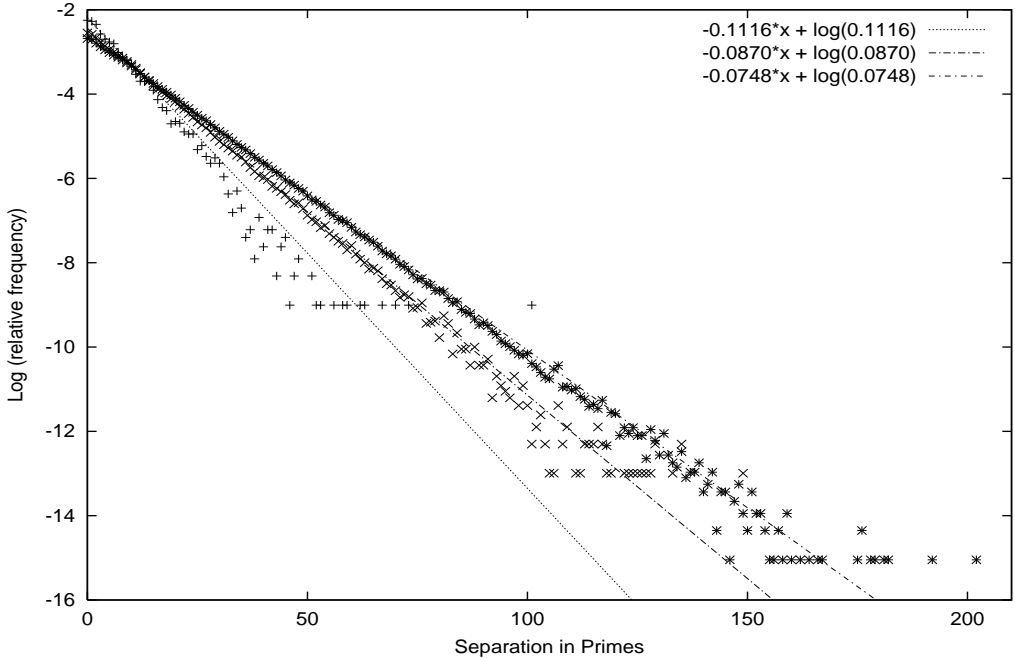}
\end{turn}
\end{center}
\caption{Data and Linear Fits for log(frequency) {\it vs.}
separation, in the cases $N = 10^6, 10^8, 10^9$.}
\label{graph2}
\end{figure}

We qualify this statement, however, since the curve in Figure~\ref{graph1} 
is only illustrative because the slope varies inversely with $N$.
Figure~\ref{graph2} displays curves (with associated best-fit straight lines) 
for the twin prime separation data for ranges $[ 5 , 10^6]$, $[ 5 , 10^8]$ 
and $[ 5 , 10^9]$, illustrating the variation with $N$. 

\noindent
Were it not the case that the magnitude of the slope diminished for 
larger values of $N$ then the Hardy--Littlewood Conjecture for twins 
would certainly fail to hold.
If the slopes of our lines were indeed the {\it same} value for 
all $N$, meaning that the probability of a given prime being a 
member of a twin pair is a universal constant, then the number 
of twin primes $\pi_2(N)$ would just be a fixed fraction of $\pi_1(N)$ 
in disagreement with Hardy--Littlewood and the empirical data.

The linear fits that we employed were constrained to ensure that the 
relative frequencies are properly normalized.
That is, if the relative frequency with which separation $s$ occurs
obeys the exponential relation consistent with our data, then for 
the sum of frequencies to be normalized to $1$, we must have
\begin{equation}
+ \mbox{(intercept)} \equiv \ln( -\mbox{(slope)})\, ,
\end{equation}
so the fit that we performed was constrained by the one-parameter 
{\it Ansatz} 
\begin{equation}
f(s) = - m s + \ln(m)\, .
\end{equation}
Table~\ref{table1} gives the best-fit values of this probability-conserving 
slope ($m$) for various $N$, along with simple statistical estimates of 
the uncertainty.
The quoted error estimates measure only the quality of the estimate 
of $m$ and and does not account for effects resulting 
from our arbitrary choices of $N$.
It may very well be the case that a more realistic assessment of
error would double the values quoted in the table.

\begin{table}[htb]
\begin{center}
\begin{tabular}{|c||c|c|c|c|}
\hline
$\pi_2(N)$ & slope & stat.~error ($\pm$) & $\pi_1(N)$ & $N$ \\
\hline
\hline
$1 \times 10^3$ & 0.141667 & 0.00599 &      7793 &      79561 \\
$5 \times 10^3$ & 0.122415 & 0.00315 &     45886 &     557521 \\
$1 \times 10^4$ & 0.114097 & 0.00325 &     97255 &    1260991 \\
$5 \times 10^4$ & 0.104126 & 0.00105 &    556396 &    8264959 \\
$1 \times 10^5$ & 0.096421 & 0.00095 &   1175775 &   18409201 \\
$5 \times 10^5$ & 0.086700 & 0.00056 &   6596231 &  115438669 \\
$1 \times 10^6$ & 0.081143 & 0.00041 &  13804822 &  252427603 \\
$3 \times 10^6$ & 0.075491 & 0.00035 &  44214960 &  863029303 \\
$5 \times 10^6$ & 0.073150 & 0.00031 &  75860671 & 1523975911 \\
$8 \times 10^6$ & 0.070965 & 0.00032 & 124538861 & 2566997821 \\
$1 \times 10^7$ & 0.070154 & 0.00029 & 157523559 & 3285916171 \\
$1.2 \times 10^7$ & 0.069814 & 0.00024 & 190894477 & 4020634603 \\
\hline
\end{tabular}
\caption{Values of slope, statistical error, $\pi_{1}(N)$, and $\pi_{2}(N)$ 
for certain $N$ from $79561$ to $4020634603$.}
\label{table1}
\end{center}
\end{table}

Two comments must be made.
The first is that all of the separations which appeared in the 
data received equal (frequency-weighted) consideration in 
our computation of best-fit slopes.
This has a consequence insofar as the large-separation, low-frequency
events constituting the tail of the distribution reduce the magnitude
of the slope, as is readily seen in Figures~\ref{graph1} and \ref{graph2}.
One might well be inclined to truncate the data by excising the
tails and fixing the slopes by the (more-strongly-linear) low-separation
data for each $N$.
We did not do this because it would have entailed a generally systematic
discarding of data from pairs appearing near the upper limit of the
range, and thus would nearly correspond to the slope with greater magnitude
that one would expect associated with an {\it effective} upper
limit $N_{{\hbox{\small eff}}} < N$.
Viewed from this perspective, it is better to consider all points
rather than submit to this degree of uncertainty.
The second comment is that we have thus far adhered to the convention
of expressing all of our results in terms of the natural number $N$.
We shall now pass over to a characterization in terms of $\pi_{1}$
-- itself a function of $N$ -- which better suits our viewpoint
that the analysis of the distribution of twins is most meaningful
when considered in terms of the primes themselves.

\begin{figure}[htb]
\begin{center}
\begin{turn}{-90}
\leavevmode
\epsfxsize=3.5in
\epsfysize=5.5in
\epsfbox{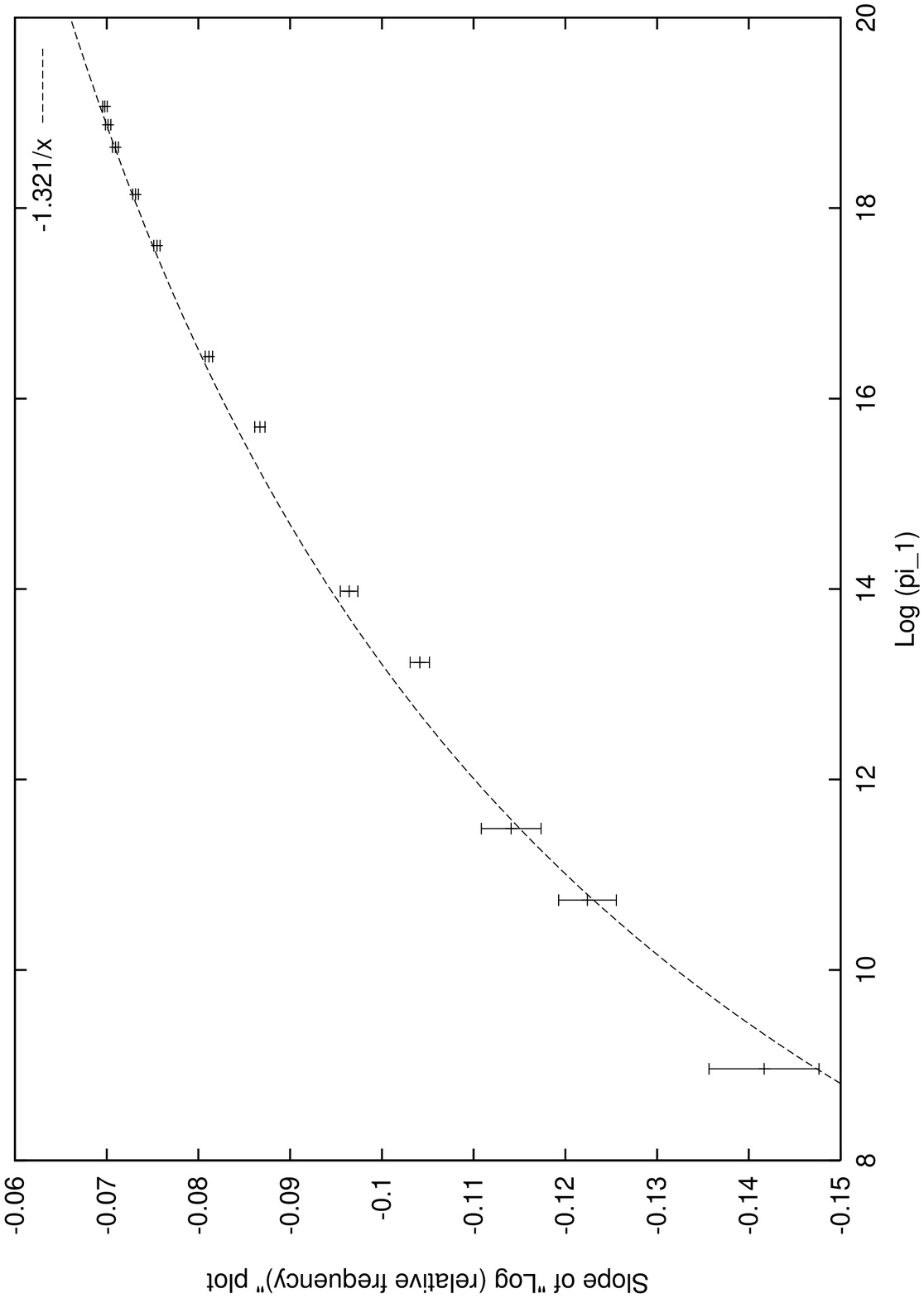}
\end{turn}
\end{center}
\caption{Slope Data and our empirical fit~(\ref{cool-eh}) {\it vs.}
$\log(\pi_{1})$.}
\label{graph3}
\end{figure}
In light of the above comments, we sketch in Figure~\ref{graph3} a plot of 
the slopes (computed in the manner described) {\it versus} 
$\log( \pi_{1} )$.
The trend seen on the graph may be well-described by
the function (remember that the error bars are understated)
\begin{equation}
- m(x) \simeq - (1.321 \pm 0.010) \big/ x \ , \quad 
\mbox{for  } x = \log( \pi_{1}(N) ) \, .
\label{cool-eh}
\end{equation}
There are two amazing features of this functional form for the dependence
of the slope on $\pi_{1}$.
The first is that the factor which appears looks suspiciously like
$-2 c_2$, the twin primes constant!
This result will be confirmed in the next section.
The second feature is that 
\[
\lim_{x \longrightarrow \infty} - m(x) = 0^{-} \ ,
\]
{\it i.e.,} as one progresses through the 
infinite set of primes, the slope which governs the distribution 
of twins does not crash through zero.
This is consistent with the Hardy--Littlewood conjecture
insofar as the twins become progressively more sparse within the
set of primes.
In addition, it is consistent with the Twin Prime Conjecture
in that the reciprocal of the slope admits the interpretation of
being the ``expected number of primes interspersed between a given twin 
and the next twin in order.''
\begin{equation}
\bar{s} = \frac{1}{m} \ , \quad \mbox{for random fixed probability events.}
\end{equation}
Since $\bar{s}$ remains finite for all $\pi_{1}$, then we can
conjecture that wherever one happens to be in the infinite set of 
primes it is possible to characterize the expected number of
primes that will be encountered on the way to the next twin.
This is also consistent with the Twin Prime Conjecture, although
unfortunately it is empirical and does not constitute a proof.

\section{Interpretive Framework}  

Recall that, up to this point, all of our results have been
purely empirical.
Now, we will argue for their essential truth and consistency
beyond the range of our data.

Consider our approximation~(\ref{pi1}) for $\pi_{1}(N)$.
Making the additional draconian approximation that the integrand
is constant at its minimum value, and discarding a small term, we get
the oft-quoted estimate 
\begin{equation}
\widetilde{\pi}_{1}(N) \sim \frac{N}{\log(N)}\, .
\label{blah1}
\end{equation}
In precisely the same manner we get
\begin{equation}
\widetilde{\pi}_2 (N) \sim 2 c_2 \, \frac{N}{\left(\log(N) \right)^2}\, .
\label{blah2}
\end{equation}

The entire set of prime numbers less than $N$ consists of the 
$2 \times \pi_{2}(N)$ elements which occur together in twins and 
$\pi_{1}(N) - 2 \times \pi_{2}(N)$ singletons.
Now, let us suppose that the set of twins is randomly interspersed 
among the singletons.
This would imply that between each twin pair there will appear,
on average, $s_{0}(N)$ singletons, where
\begin{equation} 
s_{0}(N) = \frac{ \pi_{1} - 2 \pi_{2} }{\pi_{2}} \ .
\label{blah3}
\end{equation}
Note that there is an essential distinction between $\bar{s}$ which 
arises from the actual distribution of separations and $s_0$ 
which, in effect, assumes that the twins are evenly spaced.
That is, from a value of $\bar{s}( \pi_{1}(N) )$ one can infer
$m$ and thus the probability distribution of twin prime separations
characteristic of the set of primes less than $N$.
On the other hand, $s_{0}(N)$ is an average value in which no account
is taken of the details of the distribution and thus no more information
is contained in it.

In the approximation scheme developed in (\ref{blah1}) and (\ref{blah2}),
\begin{equation}
\widetilde{s}_{0} = \frac{\log(N) - 4 c_2}{2 c_2}
\simeq \frac{\log(N)}{2 c_2} \, .
\label{10}
\end{equation}
Furthermore, to very lowest-order
\begin{equation}
\log(N) \simeq \log( \widetilde{\pi}_{1} ) 
+ \log( \log( \widetilde{\pi}_{1} )) \, ,
\label{11}
\end{equation}
and hence
\begin{equation}
\widetilde{s}_{0} \approx \frac{ \log( \widetilde{\pi}_{1} )}{2 c_2} \, .
\label{12}
\end{equation}
Finally, taking this as the expected number of singleton primes 
occurring between twin prime pairs for numbers less than or equal to $N$
we see immediately that
\begin{equation}
\widetilde{m} = \frac{1}{\widetilde{s}_{0}} \approx
\frac{2 c_2}{\log(\widetilde{\pi}_{1})}
\label{13}
\end{equation}
is completely consistent with the Prime Number Theorem, 
the Hardy--Littlewood Conjecture, and with our empirical results.

As an aside, one might consider the effect of attempting to improve 
upon the draconian approximation.
It turns out that any reasonable improvement merely results in the addition
of (small) constant terms which may be neglected in the limit of large $N$.

We are quite surprised that our empirical results yield the large $N$
limit with such accuracy.

Another test of the general consistency of our model is by comparison
with $m_0$, where
\begin{equation}
m_{0}(N) = \frac{1}{s_{0}(N)} = 
\frac{\pi_{2}(N)}{\pi_{1}(N) - 2 \pi_{2}(N)} \, .
\label{14}
\end{equation}
Bearing in mind that we are modeling more accurately the actual 
distribution of prime separations with $\bar{s}$ than with $s_0$, 
we do not expect perfect agreement, but rather that the 
general trend exposed by $s_0$ will be followed by $\bar{s}$ if
investigated beyond the range thus far examined.
We sketch below a plot of $m_{0}$ {\it vs.} $\log{\pi_{1}}$ 
using precise values for $\pi_1$ and $\pi_2$ computed by 
T.R.~Nicely~\cite{nicely2}.
Note that we have made the small adjustments of decrementing 
the published $\pi_{2}$'s by one 
and decrementing the $\pi_{1}$ by two 
to take into account our skipping the anomalous prime $2$ 
and twin $(3\ 5)$.
We are quite encouraged by the correspondence of the data on 
graph, and believe that the distributional model does extend itself 
well beyond the range of our present data.

\begin{figure}[htb]
\begin{center}
\begin{turn}{-90}
\leavevmode
\epsfxsize=3.5in
\epsfysize=5.5in
\epsfbox{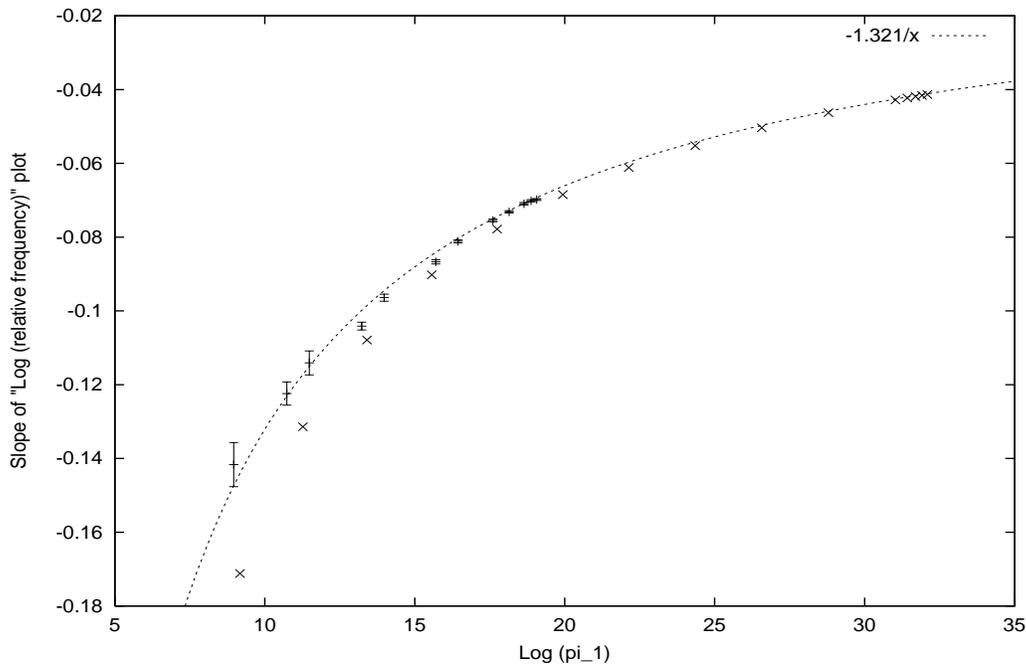}
\end{turn}
\end{center}
\caption{$m_{0}$ using Nicely Data and our empirical fit~(\ref{cool-eh}) 
{\it vs.} $\log(\pi_{1})$.}
\label{graph4}
\end{figure}

\section{Conclusion}

We believe that we have constructed a novel characterization of the 
distribution of twin primes.
The most essential feature of our approach is that we consider the 
spacings of twins among the {\it primes} themselves, rather than
among the natural numbers.
Secondly, we modeled the distribution empirically -- 
without preconceptions -- and argued that for any given $N$ 
(larger than $10^4$, say) the twin primes appear amongst the sequence 
of primes in a manner characteristic
of a completely random, fixed probability system.
Again working empirically, we noted that the ``fixed'' probability
varied with $N$, in a manner consistent with Theorems and with
Conjectures that are believed to hold.
We have parameterized the variation of the ``separation constant''
in terms of $\pi_1$, as suggested by our outlook, and have discovered 
that it has a particularly simple functional form and is also
consistent with the established Theorems and Conjectures.

With this model for the distribution now in hand and assumed viable, 
we are beginning to investigate other consequences.
These will be reported upon in a forthcoming paper~\cite{pktpinprep}.

\section{Acknowledgements}

PFK and TP thank J.~Calvo and J.~Coykendall for commenting upon a draft 
version of this manuscript.
This work was supported in part by the National Science Foundation
(USA) under grant \#OSR-9452892 and a Doctoral Dissertation 
Fellowship.

\end{document}